\documentclass[12pt]{amsart}
\usepackage{latexsym, amsmath, amscd, amssymb, amsthm}

\newtheorem*{te}{Theorem}
\newtheorem{lm}{Lemma}

\begin{document}

\noindent

 \title{  On complete system of covariants for the   binary form of  degree 8}

\author{Leonid Bedratyuk} \address{ Khmelnytsky National University, Instytuts'ka st. 11, Khmelnytsky , 29016, Ukraine}
\email {bedratyuk@ief.tup.km.ua}
\begin{abstract}
A minimal system of homogeneous generating elements of the algebra of
covariants for the binary form of degree 8 is calculated.
\end{abstract}

\maketitle

\section{introduction}

Let $V_n$ be a vector $\mathbb{C}$-space of the binary forms of degree $d$ considered with natural action of the group  $G=SL(2,\mathbb{C}).$ Let us extend the action of the group  $G$ to the polynomial functions algebra  $\mathbb{C}[V_d \oplus \mathbb{C}^2 ].$ 
Denote  by  ${C_d=\mathbb{C}[V_d \oplus \mathbb{C}^2 ]^{\,G}}$ the corresponding subalgebra of   $G$-invariant functions. 
In the vocabulary of classical invariant theory  the algebra    $C_d$ is called     the  algebra covariants of the  binary form of d-th degree.
Let  $C_d^{+}$ be an ideal of  $C_d$ generated by all homogeneus elements of positive power. Denote by  $\bar C_d$ a set of  homogeneus elements of $C_d^{+}$ such that their images in  $C_d^{+}/(C_d^{+})^2$ form a basis of the vector space. The  set  $\bar C_d$ is caled complete system of covariants of the  d-th degree binary form. Elements of   $\bar C_d$ form a minimal system of homogeneous generating elements of the invariants algebra $C_d.$ Denote by $c_d$ a number of  elements of the set    $\bar C_d.$

The complete systems of  covariants was a topic of major research interest in classical invariant theory of the 19th century.
It is easy to show that $c_1=0, $  $c_2=2,$ $c_3=4.$ A complete system of covariants  in the case   $d=4$ was calculated by Bool, Cayley, Eisenstein, see survey \cite{Dix}.
The complete systems  of invariants and covariants in the cases   $d=5,6$ were calculated by Gordan, see  \cite{Gor}. In particular, $c_4=5,$  $c_5=23,$ $c_6=26.$

Gall's attempt  \cite{Gall-1}  to discover the complete system for the 
 case $d=7$ was unsuccessful. He did offered a system of 151 covariants but 
 the system was not a minimal system, see \cite{DL},\cite{Aut-1}. Also, 
 Sylvester's attempts \cite{SF}, \cite{Sylv-10} to find the cardinality 
 $c_7$ and covariant's degree-order distribution  of the complete system 
 were mistaken, see \cite{DL}.  Therefore, the problem of finding a minimal system of 
 homogeneous generating elements (or even a cardinality of the system) of 
 the  algebra  of covariants for  the  binary form of degree 7 is still 
 open.

 The case  $d=8$  was studied by Sylvester and Gall but they have obtained 
 completelly different results. By using Sylvester-Cayley technique, 
 Sylvester in   \cite{SF} got that  ${c_8=69.}$   Gall in    \cite{Gall}, 
 evolving the Gordan's constructibe method, offered   68 covariants as a 
 complete system of covariants for  the case  $d=8.$  Also they have got a 
 different  degree-order distribution of the covariants, see in 
 \cite{Sylv-G}, a part of their long  discussion.

 For the cases $d=9,10$ Sylvester in \cite{Sylv-10} calculated the 
 cardinalities  $c_9$ and  $c_{10}$  but the present author, using a 
 computer, found numerous mistakes in those computations.

 Therefore, the 
 complete systems of covariants for the binary form are so far 
 known only up to degree 6, see \cite{Olver}, \cite{Dix}.

 For the case $d=8$ Gall have found  the 68 invariants  in implicit way by the symbolic method. A covariant has very simple symbolic representation but it is too hard to check out and be sure if the covariant is an irreducible one. The verification could be done for an covariants in their explicit representation but some of offered by Gall invariants is impossyble even for  computer calculation due to hight tranvectant's order.

To solve the computation problem we offer a form of represenation of the covariants, which  is an intermediate form between the highly unwieldy explicit represenation and too "compressed" symbolic represenation. Also we have found new symbolic representation of the covariants, which is different from Gall's representation. The representation    use  transvectants of  low orders. A first step in the simplifycation of a calculation is  calculation a semi-invariants instead calculation of covariants. Let us consider a covariant  as a polynomial  of generating  functions of the polynomial functions algebra  $\mathbb{C}[V_d \oplus \mathbb{C}^2 ].$ Then a semi-invariat is just a leading coeficient of the polinomial with respect to usual lexicographical ordering, see \cite{Gle}, \cite{Olver}. A semi-invariant is an invariant of upper unipotent matrix  subalgebra of Lie algebra $\mathfrak{sl_{2}}.$

Let us indentify the algebra $\mathbb{C}[V_d]$  with the algebra $\mathbb{C}[X_d]:=\mathbb{C}[t,x_1,x_2,\ldots, x_d],$ and the algebra $\mathbb{C}[V_d \oplus \mathbb{C}^2 ]$ identify  with the polynomial algebra $\mathbb{C}[t,x_1,x_2,\ldots, x_n,Y_1,Y_2].$

The generating elements  $\Bigl( \begin{array}{ll}  0\, 1 \\ 0\,0 \end{array} \Bigr),$ $\Bigl( \begin{array}{ll}  0\, 0 \\ 1\,0 \end{array} \Bigr)$ of the tangent Lie algebra  $\mathfrak{sl_{2}}$ act    on  $\mathbb{C}[V_d]$  by  derivations
$$
\begin{array}{l}
D_1:=\displaystyle t\frac{\partial}{\partial x_1}+2\, x_1\frac{\partial}{\partial x_2}+\cdots +d\,x_{d-1}\frac{\partial}{\partial x_d}, \\
D_2:=\displaystyle d\,x_1\frac{\partial}{\partial t}+(d-1)\,x_2\frac{\partial}{\partial x_1}+\cdots +x_{d}\frac{\partial}{\partial x_{d-1}}.
\end{array}
$$

It follows that the semi-invariant algebra   coincides with an algebra $\mathbb{C}[X_d]^{D_1},$ of polynomial solutions of the following first order PDE, see \cite{Hilb}, \cite{Gle}:
$$
\begin{array}{r}
\displaystyle t\frac{\partial u}{\partial x_1}+2\, x_1\frac{\partial u}{\partial x_2}+\cdots +d\,x_{d-1}\frac{\partial u}{\partial x_d}=0, \\
\end{array}
\eqno (*)
$$
 where $u \in \mathbb{C}[X_d],$  and  $$\mathbb{C}[X_d]^{D_1}:=\{ f \in \mathbb{C}[X_d]| D_1(f)=0 \}.$$

 It is easy to get an explicit form of the algebra  $\mathbb{C}[X_d]^{D_1},$ see, for example, \cite{Ess}. Namely -- $$\mathbb{C}[X_d]^{D_1}=\mathbb{C}[t, z_2,\ldots, z_d][\frac{1}{t}] \cap \mathbb{C}[X],$$ here  $z_i$ are some  functional independed semi-invariants of power  $i.$ The polynomials $z_i$ arised in the first time in Cayley, see  \cite{Gle}. Therefore,  any semi-invariant we may write as rational fraction of ${\displaystyle \mathbb{C}[Z_d][\frac{1}{t}]:=\mathbb{C}[t, z_2,\ldots, z_d][\frac{1}{t}].}$  This form of semi-invariants is more compact than their standard form as a polynomial of   $\mathbb{C}[X_d].$ In this case, a semi-invariant has terms number  in tens times less than terms number of the  corresponding covariant that make a computation crucial easy. 
By using Robert's theorem, see \cite{Rob}, knowing a semi-invariant one may restore a corresponding  covariant.

 In the paper on explicit form have found a complete system of covariants 
 for the  binary form of  degree 8. The system consists of 69 covariants, 
 i.e. $c_8=69.$ Also,  the covariants degree-order distribution coincides 
 completelly with Sylvester's offered   distribution.

All calculation were done with Maple.
\section{Premilinaries. }


Before any calculation we try make a simplification of a covariant represenation and their computation. 
Let  ${\it \kappa} : C_d \longrightarrow \mathbb{C}[X_d]^{D_1}$ be the   $\mathbb{C}$-linear map takes each homogeneous covariant of order $k$ to his leading coefficient, i.e. a coefficient of   $Y_1^k.$ Follow by classical  tradition an element of the algebra   $\mathbb{C}[X_d]^{D_1}$  is called  {\it semi-invariants, }  a degree of a homogeneous covariant with respect to the variables set $X_d$ is called {\it degree} of the covariant and its  degree with respect to the variables set $Y_1, Y_2$  is called {\it order.}

Suppose  $F=\displaystyle \sum_{i=0}^m \, f_i { m \choose i } Y_1^{m{-}i} Y_2^i$ be a covariant of order $m,$  $\kappa(F)=f_0 \in \mathbb{C}[X_d]^{D_1}.$ The classical Robert's theorem, \cite{Rob},  states that the covariant $F$ is completelly and   uniquely determined  by its leading coefficient $f_0,$ namely 
$$
F=\sum_{i=0}^{m} \frac{D_2^i(f_0)}{i!} Y_1^{m-i}Y_2^i.
$$
On the other hand, every semi-invariant is a leading coefficient of some covariant, see \cite{Gle}, \cite{Olver}. This 
give us well defined   explicit form of the inverse map  $${\kappa^{-1} :   \mathbb{C}[X]^{d_1} \longrightarrow C_d,}$$ namely

$$
\kappa^{-1}(a)=\sum_{i=0}^{{\rm ord}(a)} \frac{D_2^i(a)}{i!} Y_1^{{\rm ord}(a)-i}Y_2^i,  
$$
here  $ a \in  \mathbb{C}[X]^{d_1}$ and  ${\rm ord}(a)$ is an order of the element $a$ with respect to the locally nilpotent derivation  $D_2,$  i.e. ${\rm ord}(a):=\max \{ s, D_2^s(a) \neq 0 \}.$  For example, since ${\rm ord}(t)=d,$ we have 
$$
\kappa^{-1}(t)=\sum_{i=0}^{{\rm ord}(t)} \frac{D_2^i(t)}{i!} Y_1^{{\rm ord}(t)-i}Y_2^i =t Y_1^d+\sum_{i=1}^{d} { d \choose i } x_i Y_1^{d-i}Y_2^i.
$$
As we see  the $\kappa^{-1}(t)$ is just the  basic binary form. From polynomial functions point of view the covariant $\kappa^{-1}(t)$ is the evaluation  map.

Thus,  the  problem of finding of complete system of the  algebra   $\overline C_d$ is equivalent to the  problem  of finding of complete system of semi-covariants's algebra   $\mathbb{C}[X]^{D_1}.$  It is well known classical results.

A structure of  constants algebras for such locally nilpotent  derivations can be easy determined, see for example  \cite{Ess}.  In particular, for the derivation $D_1$  we get 
$$
\mathbb{C}[X_d]^{D_1}=\mathbb{C}[t,\sigma(x_2),\dots ,\sigma(x_d)][\frac{1}{t}] \cap \mathbb{C}[X_d],
$$
\noindent
where   $\sigma: \mathbb{C}[X_d] \to \mathbb{C}(X_d)^{D_1} $  is a ring homomorphism defined by  
$$
\sigma(a)=\sum_{i=0}^{\infty} d_1^{\,i}(a) \frac{\lambda^i}{i!}, \lambda = -\frac{x_1}{t}.
$$
After not complicated simplification we obtain  $\displaystyle \sigma(x_i)=\frac{z_{i+1}}{t^i},$ where  $z_i \in  \mathbb{C}(X_d)^{D_1}$ and 
$$
z_i:= \sum_{k=0}^{i-2} (-1)^k {i \choose k} x_{i-k}  x_1^k t^{i-k-1} +(i-1)(-1)^{i+1} x_1^i, i=2,\ldots,d.
$$
Especially
$$
\begin{array}{l}
z_2={x_{2}}\,t - {x_{1}}^{2}
\\
z_3={x_{3}}\,t^{2} + 2\,{x_{1}}^{3} - 3\,{x_{1}}\,{x_{2}}\,t
\\
z_4={x_{4}}\,t^{3} - 3\,{x_{1}}^{4} + 6\,{x_{1}}^{2}\,{x_{2}}\,t - 4
\,{x_{1}}\,{x_{3}}\,t^{2}
\\
z_5={x_{5}}\,t^{4} + 4\,{x_{1}}^{5} - 10\,{x_{1}}^{3}\,{x_{2}}\,t + 
10\,{x_{1}}^{2}\,{x_{3}}\,t^{2} - 5\,{x_{1}}\,{x_{4}}\,t^{3}
\\
z_6={x_{6}}\,t^{5} - 5\,{x_{1}}^{6} + 15\,{x_{1}}^{4}\,{x_{2}}\,t - 
20\,{x_{1}}^{3}\,{x_{3}}\,t^{2} + 15\,{x_{1}}^{2}\,{x_{4}}\,t^{3}
 - 6\,{x_{1}}\,{x_{5}}\,t^{4}
\\
z_7={x_{7}}\,t^{6} + 6\,{x_{1}}^{7} - 21\,{x_{1}}^{5}\,{x_{2}}\,t + 
35\,{x_{1}}^{4}\,{x_{3}}\,t^{2} - 35\,{x_{1}}^{3}\,{x_{4}}\,t^{3}
 + 21\,{x_{1}}^{2}\,{x_{5}}\,t^{4} - 7\,{x_{1}}\,{x_{6}}\,t^{5}\\

{z_{8}}=28\,{x_{1}}^{6}\,{x_{2}}\,t {-} 56\,{x_{1}}^{5}\,{x_{3}}\,t
^{2} - 56\,{x_{1}}^{3}\,{x_{5}}\,t^{4} + 28\,{x_{1}}^{2}\,{x_{6}}
\,t^{5} - 8\,{x_{1}}\,{x_{7}}\,t^{6} {-} 7\,{x_{1}}^{8} + 70\,{x_{1
}}^{4}\,{x_{4}}\,t^{3}
\mbox{} + {x_{8}}\,t^{7} 

\end{array}
$$
Thus we obtain   
$$
\mathbb{C}[X_d]^{D_1}=\mathbb{C}[t,z_2,\ldots,z_d][\frac{1}{t}] \cap \mathbb{C}[X_d].
$$
Hence, a generating elements of the semi-invariant algebra  $\mathbb{C}[X_d]^{D_1}$ we may looking as a rational fraction $\displaystyle \frac{f(z_2,\ldots,z_d)}{t^s},$ $f \in \mathbb{C}[Z_d]:=\mathbb{C}[t,z_2,\ldots,z_d],$ $s \in  \mathbb{Z_{+}}.$

To make of a calculation with an invariants in such representation we need know an  action of the operator  $D_2$ in new coordinates $t,z_2,\ldots,z_d.$
  Denote by  $D$  extention  of the derivation  $D_2$  to the algebra $\displaystyle \mathbb{C}[Z_d][\frac{1}{t}]:$
$$
D:=D_2(t)\,\frac{\partial}{\partial t}+D_2(z_2)\,\frac{\partial}{\partial z_2}+\ldots +D_2(z_d)\,\frac{\partial}{\partial z_d}.
$$ 
In  \cite{Aut} the autor proved  that 
$$
\begin{array}{l}
D(t)=-n t \lambda,\\

D(\sigma(x_2))=\displaystyle (n-2)\sigma(x_3)-(n-4)\sigma(x_2) \lambda,\\
D(\sigma(x_i))=\displaystyle (n-i)\sigma(x_{i+1})-(n-2i)\sigma(x_i) \lambda-i(n-1) \frac{\sigma(x_2) \sigma(x_{i-1})}{t}, \mbox{ for }  i>2.
\end{array}
$$
 Taking into account   $\displaystyle \sigma(x_i)=\frac{z_{i+1}}{t^i},$ $ \displaystyle \lambda = -\frac{x_1}{t}$  we can obtain an expression and for  $D(z_i),$ ${ i=2,\ldots,d.}$ Especially, for  $d=8$ we get :
$$ 
\begin{array}{l}

D=\displaystyle 7\,{x_{1}}\,{\frac {\partial }{\partial t}} - 
\displaystyle \frac {( - 15\,{x_{1}}\,{z_{3}} + 18\,{z_{2}}^{2}
 - 4\,{z_{4}})}{t}{\frac {\partial }{\partial {z_{3}}}} 
 + {\displaystyle \frac {(20\,{x_{1}}\,{z_{4}} - 24\,{z_{2}}\,{z
_{3}} + 3\,{z_{5}})}{
t}}\,{\frac {\partial }{\partial {z_{4}}}} +  \\
\\
\mbox{} + {\displaystyle \frac {(2\,{z_{6}} + 25\,{x_{1}}\,{z_{5}
} - \displaystyle 30\,{z_{2}}\,{z_{4}})}{t}}\,{\displaystyle \frac {\partial }{\partial {z_{5}}}}
  + {\displaystyle \frac {({z_{7}} + 30\,{x_{1}}\,{z_{6}
} - 36\,{z_{2}}\,{z_{5}})}{t}} \,{\displaystyle \frac {\partial }{\partial {z_{6}}}}+
 \\
\\
\mbox{} + {\displaystyle \frac {7\,(5\,{x_{1}}\,{z_{7}} - 6\,{z_{
2}}\,{z_{6}})}{t}} \,{\displaystyle \frac {\partial }{\partial {z_{7}}}}
 + {\displaystyle \frac {5\,(2\,{x_{1}}\,{z_{2}} + {z_{3}})}{t}} \,
{\displaystyle \frac {\partial }{ \partial {z_{2}}}} 
- {\displaystyle \frac {8\,( - 6\,{x_{1}}\,{z_{8}} + 7\,{z_{2}}
\,{z_{7}})\,}{t}} {\displaystyle \frac {\partial }{ \partial {z_{8}}}}.
\end{array}
$$

To calculate the semi-invariants we need have an analogue of the transvectants. 
Suppose $$F=\sum_{i=0}^m \, f_i { m \choose i } Y_1^{m{-}i} Y_2^i, \mbox{   }  G=\sum_{i=0}^k \, f_i { k \choose i } Y_1^{k{-}i} Y_2^i,   \mbox{   }        f _i, g_i \in \mathbb{C}[Z_d][\frac{1}{t}], $$ are two covariants of the degrees  $m$ and $ k$ respectively. Let 
$$
(F,G)^r=\sum_{i=0}^r (-1)^i { r \choose i } \frac{\partial^r F}{\partial Y_1^{r-i} \partial Y_2^i}   \frac{\partial^r G}{\partial Y_1^{i} \partial Y_2^{r-i}},
$$
be their   $r$-th transvectant. 
The following lemma give us rule how to find the semi-invariant $\kappa( (F,G)^r)$ without of direct  computing of the covariant  $(F,G)^r.$
\begin{lm}
The leading coefficient  $\kappa ((F,G)^r)$ of the covariant $(F,G)^r,  0 \leq r \leq \min(m,k)$ is calculating by the formula
$$
\kappa((F,G)^r)=\sum_{i=0}^r (-1)^i { r \choose i } \frac{D^i(\kappa(F))}{[m]_i} {\Bigl | _{x_1=0,\ldots,x_r=0}}  \frac{D^{r-i}(\kappa(G))}{[k]_{r-i}}{\Bigl | _{x_1=0,\ldots,x_r=0}},
$$ 
here  $[a]_i:=a (a-1) \ldots (a-(i-1)), a \in \mathbb{Z}.$
\end{lm}
The proof is in  \cite{Aut-1}.

Let   $f,g$ be two semi-invariants. Their numers are a polynomials of  $z_2,\ldots,z_d$ with rational coefficients. Then the semi-invariant  $\kappa((\kappa^{-1}(f),\kappa^{-1}(g))^i)$  be a fraction and their numer be a polynomial  of  $z_2,\ldots,z_d$ with rational coefficients too.  Therefore we may multiply $\kappa((\kappa^{-1}(f),\kappa^{-1}(g))^r)$ by some rational number  $q_r(f,g) \in \mathbb{Q}$ such that the numer of the expression  $ q_r(f,g) \kappa((\kappa^{-1}(f),\kappa^{-1}(g))^r)$  be now a polynomial with an integer coprime coefficients. Put 
$$
[f,g]^r:=q_r(f,g) \kappa((\kappa^{-1}(f),\kappa^{-1}(g))^r), 0 \leq r\leq \min({\rm ord}(f),{\rm ord}(g)). 
$$
The expression $[f,g]^r$ is said to be the  r-th {\it semitranvectant} of the semi-invariants $f$ and $g.$

The following statements is direct consequences of corresponding  transvectant properties,  see \cite{Gle}:
\begin{lm} Let  $f, g$ be two semi-invariants.  Then the folloving conditions hold
\begin{enumerate}
\item[({\it i})] the semitransvectant $[t, f\,g]^i$ is reducible for  $ 0 \leq i \leq \min(d, \max( {\rm ord}(f),{\rm ord}(g));$
\item[({\it ii})] if  $ {\rm ord}(f)=0,$ then   $[t, f\,g]^i=f [t,g]^i;$
\item[({\it iii})] ${\rm ord}([f,g]^i)={\rm ord}(f)+{\rm ord}(g)-2\, i;$ 
\item[({\it iv})] ${\rm ord}(z_2^{i_1}z_3^{i_3} \cdots z_d^{i_d})=d\,(i_2+i_3+\cdots +i_d)-2\,(2\, i_2+3\, i_3+\cdots +d \, i_d).$
\end{enumerate}
\end{lm}

\section{Calculation }


Let  $\overline C_{8,\,i}:= (C_8)_i$ be a subset  of  $\overline C_{8}$ whose elements has degree  $i.$ 
Let   $C_{+}$ be the  ideal $\sum_{i>0}C_{8,i}$ of   $C_8,$ and  $\overline C_{8,\,i}:=\overline C_8 \cap C_{8,i}.$
The number  $\delta_i$  of linearly independent irreducible invariants of degree $i$  is calculated by the formula $\delta_i =\dim C_{8,i}-\dim(C^2_{+})_i.$ A dimension of the vector space  $C_{d\!,i}$  is calculated  by  Sylvester-Cayley formula, see for example \cite{Hilb}, \cite{SP}:
$$
\dim C_{d,\,i}=\frac{(1-T^{d+1})(1-T^{d+2})\ldots (1-T^{d+i})}{(1-T^2)\ldots(1-T^i)}{\Bigl | _{T=1}}.
$$ 
A dimension of the vector space $(C^2_{+})_i$ is calculating by the formula 
$\dim(C^2_{+})_i=\sigma_i -\dim S_i.$ Here  $\sigma_i$ is a coefficient of  $T^i$ in Poincar\'e series  $\displaystyle \frac{1}{\prod_{k<i} (1-T^k)^{\delta_k} }$ of a graded algebra generated by the system of homogeneous elements  $\displaystyle  \cup_{k<i} \overline C_{8,k}, $ and  $S_i$ is a vector subspace of  $(C^2_{+})_i$ generated by syzygies. The dimension  $\dim S_i$ one  may find by direct Maple calculation.

If we already have calculated the set $\overline C_{8,\,i}$, then the  elements of the set  $\overline C_{8,\,i+1} $ we are seeking  as an irreducible elements of a  basis of a vector space  generated by semitransvectants of the  form  $[t,u\,v]^r, u \in C_{8,\,l},$  $v \in C_{8,\,k},$ $l+k=i,$ $ \max({\rm ord}(u),{\rm ord}(v)) \leq r \leq 8.$ It is a standard linear algebra problem.

The unique semi-invariant of the degree one obviously is  $t, {\rm ord}(t)=8.$

For $i=2$ have  $\dim C_{8\!,2}=5,$ $\sigma_2=1, $  thus  $\delta_{8,2}=4.$ The semi-transvectants  $[t,t]^k$  are equal to zero for odd  $i.$ Put 
$$
\begin{array}{l}
dv_1:=[t,t]^2=z_2= \,{x_{2}}\,t - {x_{1}}^{2}, 
  \\ 
dv_2:=[t,t]^4={\displaystyle \frac {{z_{4}} + 3\,{z_{2}}^{2}}{t^{2}}}={x_{4
}}\,t - 4\,{x_{1}}\,{x_{3}} + 3\,{x_{2}}^{2},  \\

dv_3:=[t,t]^6={\displaystyle \frac {{z_{6}} + 15\,{z_{2}}\,{z_{4}} - 10\,{z_{3}
}^{2}}{t^{4}}} = \,{x_{6}}\,t - 6\,{x_{1}}\,{x_{5}} + 15\,{x_{2}}
\,{x_{4}} - 10\,{x_{3}}^{2}, \\

dv_4:=[t,t]^8={\displaystyle \frac {{z_{8}} + 28\,{z_{2}}\,{z_{6}} - 56\,{z_{3}
}\,{z_{5}} + 35\,{z_{4}}^{2}}{t^{6}}} = - 8\,{x_{1}}\,{x_{7}}
 + {x_{8}}\,t + 28\,{x_{2}}\,{x_{6}} - 56\,{x_{3}}\,{x_{5}} + 35
\,{x_{4}}^{2},\\
{\rm ord}(dv_1)=12, \mbox {    } {\rm ord}(dv_2)=8, \mbox {    }  {\rm ord}(dv_3)=4,  \mbox {    }  {\rm ord}(dv_4)=0.
\end{array}
$$
The polynomials     $t^2,$ $dv_1,$ $dv_2,$ $dv_3$ $dv_4$   are linear independent. Therefore,   the set   $\overline C_{8,\,2}, $ consists of the irreducible semi-invariants  $dv_1,$ $dv_2,$ $dv_3,$ $dv_4.$

For $i=3$ we have  $\dim C_{8\!,3}=13,$ $\sigma_3=5, $ $\dim S_3=0,$ thus $\delta_{8,3}=8.$

The set $\overline C_{8,\,3}$  consists of the following 8 irreducible semi-invariants:  
$$
\begin{array}{ll}
tr_1=[t,dv_1]^3, {\rm ord}(tr_1)=14, & tr_2=[t,dv_1]^4,  {\rm ord}(tr_2)=12,\\

tr_3=[t,dv_1]^5,  {\rm ord}(tr_3)=10, & tr_4=[t,dv_1]^7,  {\rm ord}(tr_4)=6,\\

tr_5=[t,dv_2]^4,  {\rm ord}(tr_5)=8, & tr_6=[t,dv_2]^6,  {\rm ord}(tr_6)=4,\\

tr_7=[t,dv_2]^8,  {\rm ord}(tr_5)=0, & tr_8=[t,dv_1],  {\rm ord}(tr_6)=18.
\end{array}
$$
  
For $i=4$ we have  $\dim C_{8\!,4}=33,$ $\sigma_4=23, $ $\dim S_4=0,$ thus $\delta_{8,4}=10.$        
The set $\overline C_{8,\,4}$ consists of the following 10 irreducible semi-invariants:  
$$
\begin{array}{ll}
ch_1=[t,tr_4]^2, {\rm ord}(ch_1)=10, & ch_2=[t,tr_4]^5,  {\rm ord}(ch_2)=4,\\

ch_3=[t,tr_5]^8,  {\rm ord}(ch_3)=0, & ch_4=[t,tr_6]^4,  {\rm ord}(ch_4)=4,\\

ch_5=[t,tr_1]^4,  {\rm ord}(tr_5)=14, & ch_6=[t,tr_1]^5,  {\rm ord}(tr_6)=12,\\

ch_7=[t,tr_1]^6,  {\rm ord}(ch_7)=10, & ch_8=[t,tr_1]^7,  {\rm ord}(ch_6)=8,\\

ch_9=[t,tr_2],  {\rm ord}(ch_7)=18, & ch_{10}=[t,tr_2]^7,  {\rm ord}(ch_6)=6.
\end{array}
$$

For $i=5$ we have  $\dim C_{8\!,5}=73,$ $\sigma_5=65. $ 
Below is typical instance how $\dim S_i$ is calculated. 
The vector space  $(C^2_{+})_5=t \, C_{8,4}+C_{8,2}\,C_{8,3}$  is generated by the following 65 elements: 
$$
\begin{array}{l}
t\,ch_1,\ldots, t\,ch_{10},\\
t^2\,tr_1,\ldots,t^2\,tr_8,\\
dv_1\,tr_1,\ldots, dv_4\,tr_8,\\
t^3\,dv_1,\ldots,t^3\,dv_4,\\
t^5.
\end{array}
$$
To find a basis  of the vector space $S_5$ of syzygies let us  equate the system :
$$
\alpha_1 t\,ch_1+\alpha_2 t\,ch_2+\cdots+\alpha_{65} t^5=0.
$$

By solving it we get
$$
\begin{array}{l}
(55\,{\mathit{tr}_{1}}\,{\mathit{dv}_{3}} - 55\,{\mathit{tr
}_{3}}\,{\mathit{dv}_{2}} + {\mathit{ch}_{7}}\,\mathit{t}
 - 12\,{\mathit{ch}_{1}}\,\mathit{t})\,{\alpha _{19}} +\\
\mbox{} + ({\displaystyle \frac {383}{5}} \,{\mathit{tr}_{4}}\,
\mathit{t}^{2} + {\mathit{ch}_{5}}\,\mathit{t} - 
{\displaystyle \frac {176}{5}} \,{\mathit{tr}_{8}}\,{\mathit{
dv}_{3}} + {\displaystyle \frac {176}{5}} \,{\mathit{tr}_{3}}
\,{\mathit{dv}_{1}})\,{\alpha _{24}} +\\
\mbox{} + (126\,{\mathit{tr}_{1}}\,{\mathit{dv}_{1}} - {
\mathit{ch}_{9}}\,\mathit{t} + {\mathit{tr}_{3}}\,\mathit{
t}^{2} - 126\,{\mathit{tr}_{8}}\,{\mathit{dv}_{2}})\,{
\alpha _{11}} =0.
\end{array}
$$
Therefore, the vector space   $S_5$ is generated by the 3 syzygies:
$$
\begin{array}{l}
 - 12\,{\mathit{ch}_{1}}\,\mathit{t} + 55\,{\mathit{tr}_{1}
}\,{\mathit{dv}_{3}} - 55\,{\mathit{tr}_{3}}\,{\mathit{dv}
_{2}} + {\mathit{ch}_{7}}\,\mathit{t}=0,\\
5\,{\mathit{ch}_{5}}\,\mathit{t} + 383\,{\mathit{tr}_{4}}\,
\mathit{t}^{2} - 176\,{\mathit{tr}_{8}}\,{\mathit{dv}_{3}}
 + 176\,{\mathit{tr}_{3}}\,{\mathit{dv}_{1}}=0,\\
 - {\mathit{ch}_{9}}\,\mathit{t} - 126\,{\mathit{tr}_{8}}\,
{\mathit{dv}_{2}} + 126\,{\mathit{tr}_{1}}\,{\mathit{dv}_{1
}} + {\mathit{tr}_{3}}\,\mathit{t}^{2}=0.
\end{array}
$$
Thus  $\dim S_5=3$ i  $\delta_{8,5}=11$ and the set  
$\overline C_{8,\,5}$ consists of  11 irreducible semi-invariants.
We searshing the semi-invariants  as semitransvectants of the form  $[t,u]^i,$  $u \in (C^2_{+})_4.$ By using Lemma 2 we have that irreducibles covariants  be only for the following values of  $u:$ 
$$
{\mathit{ch}_{1}}, \,{\mathit{ch}_{2}}, \,{\mathit{ch}_{4}},\,{\mathit{ch}_{5}}, \,{\mathit{ch}_{6}}, \,{\mathit{ch}_{
7}}, \,{\mathit{ch}_{8}}, \, \,{\mathit{
ch}_{9}}, \,{\mathit{ch}_{10}}, \,{\mathit{dv}_{3}}^{2}.
$$

We calculate the 65 semitransvetants of the forms 
 $[t, ch_i]^k,$ $k=1,\ldots, \min(8, {\rm ord}(ch_i)),$ ${i=1,2,4,\ldots, 10,}$    $[t,dv_3^2]^k,$ $k=5,\dots,8,$  and select 11 linearly independed such that ones dont belong to the vector space  $ (C^2_{+})_5:$
$$
\begin{array}{ll}
pt_1=[t,dv_3^2]^6, {\rm ord}(pt_1)=4, & pt_2=[t,dv_3^2]^7,  {\rm ord}(pt_2)=2,\\

pt_3=[t,dv_3^2]^8,  {\rm ord}(pt_3)=0, & pt_4=[t,ch_1]^2,  {\rm ord}(pt_4)=14,\\

pt_5=[t,ch_1]^4,  {\rm ord}(pt_5)=10, & pt_6=[t,ch_1]^5,  {\rm ord}(pt_6)=8,\\

pt_7=[t,ch_1]^7,  {\rm ord}(pt_7)=4, & pt_8=[t,ch_2],  {\rm ord}(pt_6)=10, \\

pt_9=[t,ch_4],  {\rm ord}(pt_9)=10, & pt_{10}=[t,ch_4]^3,  {\rm ord}(pt_{10})=6, \\

pt_{11}=[t,dv^3_2]^5,  {\rm ord}(pt_9)=6.
\end{array}
$$

Thus, $\overline C_{8,\,5}=\{pt_1,pt_2,\ldots, pt_{11}\}.$

For $i=6$ we have $\dim C_{8\!,6}=151,$ $\sigma_5=172,$ $\dim S_6=30.$
Thus $\delta_6=151-(172-30)=9,$   the set
$\overline C_{8,\,6}$ consists of the  9 irreducible semi-invariants: 
$$
\begin{array}{ll}
sh_1=[t,tr_6\,dv_3]^6, {\rm ord}(sh_1)=4, & sh_2=[t,tr_6\,dv_3]^7,  {\rm ord}(sh_2)=2,\\

sh_3=[t,tr_6\,dv_3]^8,  {\rm ord}(sh_3)=0, & sh_4=[t,pt_5]^5,  {\rm ord}(sh_4)=8,\\

sh_5=[t,pt_6]^5,  {\rm ord}(sh_5)=6, & sh_6=[t,pt_8]^6,  {\rm ord}(sh_6)=6,\\

sh_7=[t,pt_9]^4,  {\rm ord}(sh_7)=6, & sh_8=[t,pt_9]^7,  {\rm ord}(sh_6)=4, \\

sh_9=[t,pt_10]^2,  {\rm ord}(sh_9)=10. \\
\end{array}
$$

For $i=7$ we have  $\dim C_{8\!,7}=289,$ $\sigma_7=385,$ $\dim S_7=104.$
Thus  $\delta_7=8,$  the set 
$\overline C_{8,\,7}$ consists of the  8 irreducible semi-invariants:  
$$
\begin{array}{ll}
si_1=[t,ch_{10}\,dv_3]^7, {\rm ord}(si_1)=4, & si_2=[t,tr_6^2]^5,  {\rm ord}(si_2)=6,\\

si_3=[t,tr_6^2]^7,  {\rm ord}(si_3)=2, & si_4=[t,tr_6^2]^8,  {\rm ord}(si_4)=0,\\

si_5=[t,sh_9]^6,  {\rm ord}(si_5)=6, & si_6=[t,sh_9]^7,  {\rm ord}(si_6)=4,\\

si_7=[t,ch_4\,dv_3]^5,  {\rm ord}(si_7)=6, & si_8=[t,ch_{10}\,dv_3]^8,  {\rm ord}(si_6)=2. \\

\end{array}
$$      
 For $i=8$ we have  $\dim C_{8\!,8}=289,$ $\sigma_8=385,$ $\dim S_8=104.$
Thus  $\delta_8=8.$ The set 
$\overline C_{8,\,8}$ consists of the  7 irreducible semi-invariants: 
$$
\begin{array}{ll}
vi_1=[t,ch_4 \,tr_6]^7, {\rm ord}(vi_1)=2, & vi_2=[t,ch_4\,tr_6]^8,  {\rm ord}(vi_2)=0,\\

vi_3=[t,ch_2\,tr_6]^5,  {\rm ord}(vi_3)=6, & vi_4=[t,pt_{10}\,dv_3]^7,  {\rm ord}(vi_4)=4,\\

vi_5=[t,pt_{10}\,dv_3]^8,  {\rm ord}(vi_5)=2, & vi_6=[t,ch_4\, tr_6]^5,  {\rm ord}(vi_6)=6,\\

vi_7=[t,ch_4 \,tr_6]^6,  {\rm ord}(vi_7)=4.
\end{array}
$$   
  For $i=9$ we have  $\dim C_{8\!,8}=910,$ $\sigma_8=1782,$ $\dim S_8=877.$
Thus  $\delta_8=5.$ The set 
$\overline C_{8,\,9}$ consists of the  5 irreducible semi-invariants: 
$$
\begin{array}{ll}
de_1=[t,vi_2]^6, {\rm ord}(de_1)=4, & de_2=[t,vi_2]^7,  {\rm ord}(de_2)=2,\\

de_3=[t,sh_2\,dv_3]^6,  {\rm ord}(de_3)=2, & de_4=[t,pt_1 \,tr_6]^8,  {\rm ord}(de_4)=0,\\

de_5=[t,ch_4^2]^7,  {\rm ord}(de_5)=2.
\end{array}
$$   

  For $i=10$ we have  $\dim C_{8\!,10}=1514,$ $\sigma_{10}=3673,$ $\dim S_{10}=2162.$
Thus  $\delta_{10}=3.$ The set 
$\overline C_{8,\,10}$  consists of the  3 irreducible semi-invariants: 
$$
\begin{array}{ll}
des_1=[t,pt_1 ch_4]^8, {\rm ord}(des_1)=0, & des_2=[t,si_2\,dv_3]^8,  {\rm ord}(des_2)=2,\\

des_3=[t,pt_10\,ch_4]^8,  {\rm ord}(des_3)=2.

\end{array}
$$   
  For $i=11$ we have  $\dim C_{8\!,11}=2430,$ $\sigma_{11}=7355,$ $\dim S_{11}=4927.$
Thus   $\delta_{11}=2.$ The set 
$\overline C_{8,\,11}$ consists of the  2 irreducible semi-invariants: 
$$
\begin{array}{ll}
odn_1=[t,si_3\,tr_6]^6, {\rm ord}(odn_1)=2, & odn_2=[t,vi_7\,dv_3]^7,  {\rm ord}(odn_2)=2.\\
\end{array}
$$  
   For $i=12$ we have  $\dim C_{8\!,12}=3788,$ $\sigma_{12}=14520,$ $\dim S_{12}=10733.$
Thus   $\delta_{12}=1.$ The set 
$\overline C_{8,\,12}$ consists of the unique  irreducible semi-invariant: 

$$
\begin{array}{ll}
dvan=[t,vi_5\,tr_6]^6, & {\rm ord}(dvan_1)=2. 
\end{array}
$$   
It is follows from \cite{Gall-1}  that $\delta_i=0$ for $i>12.$

The cardinalities of the sets  $\overline C_{8,\,i}$ and the covariant's degree-order distributions of Sylvester's results, see \cite{SF}.

 Summarizing the above results we get
\begin{te}
The  system of the 69 covariants:
$$
\begin{array}{l}
t,\\
dv_1,dv_2,dv_3,dv_4,\\
tr_1, tr_2,tr_3,tr_4,tr_5,tr_6,tr_7, tr_8,\\
ch_1, ch_2, ch_3, ch_4, ch_5, ch_6, ch_7, ch_8,ch_9, ch_{10},\\
pt_1, pt_2, pt_3, pt_4, pt_5, pt_6, pt_7, pt_8, pt_9, pt_{10}, pt_{11},\\
sh_1,sh_2,sh_3,sh_4,sh_5, sh_6, sh_7,sh_7,sh_9,\\
si_1, si_2, sh_3,sh_4,sh_5,sh_6,sh_7,si_8, \\
vi_1, vi_2, vi_3,vi_4,vi_5,vi_6, vi_7,\\
de_1, de_2, de_3, de_4, de_5,\\
des_1, des_2, des_2,\\
odn_1, odn_2,\\
dvan.
\end{array}
$$
is a complete system of the covariants for the  binary form of  degree 8.

\end{te}

\section{Appendix}

The degree-order distribution for $\overline C_8.$
\newpage

\begin{table}
\[ \qquad \qquad \text{order} \] \[ \text{degree} \; \; \begin{tabular}{|c||c|c|c|c|c|c|c|c|c|c|c|c|c|c|c|c|c|c|c|c|c|c|}
 \hline       &  0  & 2&  4 & 6 & 8  &10 & 12 &  14 & 16   & 18    \\ 
\hline 
\hline
 1                   & {} & {} & {} & {}  & t  &  {} & {}    &  {} & {}   & {}       \\ 
\hline
 2                   & $dv_4$ & {} & $dv_3$  & {}  & $dv_2$  &  {} & $dv_1$    &  {} & {}  & {}        \\ 

 \hline 3         & $tr_7$ & {} & $tr_6$ &$tr_4$   & $tr_5$  &  $tr_3$ & $tr_2$    &  $tr_1$ & {}   &$tr_8$       \\ 

\hline
 4                   &$ch_3$ & {} & $ch_2, ch_4$ & $ch_{10}$  & $ch_8$  &  $ch_1, ch_7$  & $ch_6$    &  $ch_5$ & {}   & $ch_9$      \\ 
   
\hline
5                    & $pt_3$ &  $pt_2$ & $pt_{1}, pt_{7}$  & $pt_{10}, pt_{11}$  & $pt_6$  &  $pt_5, pt_8, pt_9$ & {}    &  $pt_4$ & {}   & {}       \\  

 \hline
6                    &  $sh_3$ & $sh_2$  &  $sh_1, sh_8$ & $sh_5, sh_6,sh_7$  & $sh_4$  &  $sh_9$ & {}    &  {} & {}   & {}       \\   

\hline
7                    & $si_4$ & $si_3, si_8$ & $si_1, si_6$ & $si_2, si_5,si_7$   & {}  &  {} & {}    &  {} & {}   & {}     \\ 
 
 \hline
8                    &  $vi_2$ & $vi_1, vi_5$ & $vi_4, vi_7$  &$vi_3, vi_6$ & {}  &  {}  &  {}    &  {} & {}   & {}       \\    
                                                                            
 \hline
 9                   & $de_4$ & $de_2, de_3, de_5$ & $de_1$ & {}  & {}  &  {} & {}    &  {} & {}   & {}      \\ 

\hline
 10                   & $des_1$ & $des_2, des_3$ & {} & {}  & {}  &  {} & {}    &  {} & {}   & {}       \\ 

\hline
 11                   & {} & $odn_1, odn_2$ & {} & {}  & {}  &  {} & {}    &  {} & {}   & {}     \\ 

\hline
 12                   & {} & $dvan$ & {} & {}  & {}  &  {} & {}    &  {} & {}   &      \\ 

\hline
 13                   & {} &{} & {} & {}  & {}  &  {} & {}    &  {} & {}   &      \\ 

\hline \end{tabular} \] 
\end{table}

\end{document}